\documentclass[letterpaper, 10 pt, conference]{IEEEconf}  

\IEEEoverridecommandlockouts                              
\usepackage{graphicx} 
\usepackage{amsmath} 
\usepackage{amssymb}  
\usepackage{xcolor}
\usepackage{array}
\usepackage{url}

\title{\LARGE \bf
  GPU-accelerated dynamic nonlinear optimization with ExaModels and MadNLP
}

\author{François Pacaud and Sungho Shin
}

\begin{document}

\maketitle
\thispagestyle{empty}
\pagestyle{empty}

\begin{abstract}
  We investigate the potential of Graphics Processing Units (GPUs)
  to solve large-scale nonlinear programs with a dynamic
  structure.
  Using ExaModels, a GPU-accelerated automatic differentiation tool,
  and the interior-point solver MadNLP, we significantly reduce
  the time to solve dynamic nonlinear optimization problems.
  The sparse linear systems formulated in the interior-point method
  is solved on the GPU using a hybrid solver combining an iterative
  method with a sparse Cholesky factorization,
  which harness the newly released NVIDIA cuDSS solver.
  Our results on the classical distillation column instance show that despite
  a significant pre-processing time, the hybrid solver allows to reduce
  the time per iteration by a factor of 25 for the largest instance.
\end{abstract}

\section{INTRODUCTION}

There is a strong interest in using high-performance computing
for accelerating the solution of dynamic nonlinear programs, as this is critical for Nonlinear Model Predictive Control (NMPC) and optimal control applications \cite{diehl2009efficient,kirches2010efficient}.
It is well known that for problems with a dynamic structure,
the Newton step is equivalent to
a linear-quadratic problem, solvable using dynamic programming
or Riccati recursions~\cite{dunn1989efficient}.
By doing so, the method exploits the dynamic structure \emph{explicitly}.
Alternatively, one can leverage the structure \emph{implicitly} inside
a sparse direct solver, in charge of finding an appropriate ordering
to reduce the fill-in in the sparse factorization~\cite{zavala2008interior}.
That kind of many-degrees-of-freedom approaches are known to scale better with the problem's size.

\subsection{Related works}
When solving generic large-scale nonlinear programs, the two bottlenecks are
the computation of the Newton step and the evaluation of the derivatives.
On the one hand, the solution of the Newton step can be accelerated
either by exploiting the structure explicitly or by using efficient
sparse linear algebra routines~\cite{wright1991partitioned,rao1998application,cervantes1998large}.
On the other hand, efficient automatic differentiation routines
have been introduced, which now evaluate the first and second-order derivatives
in a vectorized fashion for performance~\cite{andersson2019casadi}.
As a result,
the dynamic optimization problem can be solved with near real-time
performance~\cite{verschueren2022acados} when coupled with an optimization solver~\cite{ferreau2014qpoases,frasch2015parallel,frison2020hpipm}.

With their focus on embedded applications, the solvers listed
in the previous paragraph have been heavily optimized on CPU architectures,
going as far as using dedicated linear algebra routines~\cite{frison2018blasfeo}.
Aside, NVIDIA has recently released the NVIDIA Jetson GPU, developed primarily for embedded applications. 
Hence, solving MPC problems on GPU/SIMD architectures is gaining more traction~\cite{kerrigan2015computer},
with new applications in robotics and autonomous vehicles~\cite{phung2017model,yu2017efficient,rathai2020gpu}.

\subsection{Contributions}
In this article, we investigate the capability of the
modeler ExaModels and the interior-point solver MadNLP~\cite{shin2023accelerating} --- both leveraging
GPU acceleration --- to solve nonlinear programs with dynamic structure.
MadNLP implements a filter line-search interior-point method~\cite{wachter2006implementation},
which results in solving a sequence of sparse indefinite linear systems with a saddle-point
structure~\cite{benzi2005numerical}. The linear systems are increasingly ill-conditioned as we are approaching the solution, preventing a solution with Krylov-based solvers.
The alternative is to use an inertia-revealing sparse direct solver, generally
implementing the Duff-Reid factorization~\cite{duff1983multifrontal}. Unfortunately, it
is well known that such factorization is not practical on the GPU, as they rely
on expensive numerical pivoting operations for stability~\cite{tasseff2019exploring,swirydowicz2021linear}.
The usual workaround is to densify the solution of the linear systems
using a null-space method, as was investigated in our previous work
\cite{cole2023exploiting,pacaud2022condensed}. Instead, we propose to solve the linear systems with
a hybrid sparse linear solver mixing a sparse Cholesky routine with an iterative method.
We present two alternative methods for the hybrid solver. On the one hand, Lifted-KKT
\cite{shin2023accelerating} uses an equality relaxation strategy to reduce
the indefinite linear system down to a sparse positive definite matrix, factorizable
using Cholesky factorization. On the other hand, HyKKT~\cite{regev2023hykkt} uses the Golub and Greif
method~\cite{golub2003solving} (itself akin to an Augmented Lagrangian method) to solve the linear
system with an inner direct solve used in conjunction with
a conjugate gradient.
Both methods are fully implementable on the GPU and rely only on basic linear algebra
routines. We show on the classical distillation column instance~\cite{cervantes1998large}
that despite a significant
pre-processing time, both Lifted-KKT and HyKKT reduce the time per IPM iteration
by a factor of 25 and 18 respectively, compared to the HSL solvers.

\section{PROBLEM FORMULATION}

We formulate the dynamic nonlinear optimization problem
as a generic nonlinear program.
Let $n$ the number of variables.
The objective is encoded by a function $f: \mathbb{R}^n \to \mathbb{R}$,
the dynamic and algebraic constraints by a function $g: \mathbb{R}^n \to \mathbb{R}^{m_e}$
and the remaining inequality constraints by a function $h: \mathbb{R}^n \to \mathbb{R}^{m_i}$.
By introducing slack variables $s \geq 0$, the problem writes:
\begin{equation}
  \label{eq:problem}
  \begin{aligned}
    \min_{x \in \mathbb{R}^n, s\in \mathbb{R}^{m_i}} \; & f(x) \\
    \text{s.t.} ~ & g(x) = 0 \;, ~ h(x) + s = 0  \;, ~ s \geq 0 .
  \end{aligned}
\end{equation}
We note $y \in \mathbb{R}^{m_e}$ (resp. $z \in \mathbb{R}^{m_i}$) the multiplier
attached to the equality constraints (resp. the inequality
constraints).
The Lagrangian of \eqref{eq:problem} is defined as
\begin{equation}
  L(x, s, y, z, \nu) = f(x) + y^\top g(x) + z^\top (h(x) + s)
  - \nu^\top s \; .
\end{equation}
A primal-dual variable $w := (x, s, y, z, \nu)$ is
solution of Problem~\eqref{eq:problem} if it satisfies
the Karush-Kuhn-Tucker (KKT) equations
\begin{equation}
  \label{eq:kkt}
  \left\{
  \begin{aligned}
    & \nabla_x L(x, s, y, z, \nu) = 0 \;, \\
    & z - \nu = 0 \\
    & g(x) = 0\;, \\
    & h(x) + s = 0 \;, \\
    & 0 \leq s \perp \nu \geq 0 \;,
    \end{aligned}
  \right.
\end{equation}
where we use the symbol $\perp$ to denote the complementarity
constraints $s_i \times \nu_i = 0$ for $i = 1, \cdots, m_i$.

We note the active set $\mathcal{B}(x) = \{ i = 1,\cdots, m_i \; | \; h_i(x) = 0 \}$,
and denote the active Jacobian as $A(x) = \big[ \nabla_x g(x)^\top ~ \nabla_x h_{\mathcal{B}}(x)^\top \big]^\top$
with $h_{\mathcal{B}}(x) := \{h_i(x)\}_{i \in \mathcal{B}(x)}$.
We suppose the following assumption holds
at a primal-dual solution $w^\star := (x^\star, s^\star, y^\star, z^\star, \nu^\star)$ of \eqref{eq:problem}.
\begin{itemize}
  \item Linear Independance Constraint Qualification (LICQ):
    the active Jacobian $A(x^\star)$ is full row-rank.
  \item Strict complementarity (SCS): for every $i \in \mathcal{B}(x^\star)$,
    $z_i^\star > 0$.
  \item Second-order sufficiency (SOSC): for every
    $v \in \text{null}(A(x^\star))$, $v^\top (\nabla^2_{x x} L(w^\star)) v > 0$.
\end{itemize}

\section{INTERIOR-POINT METHOD}
The primal-dual interior-point method (IPM) reformulates
the non-smooth KKT conditions~\eqref{eq:kkt} using an
homotopy method \cite[Chapter 19]{nocedal_numerical_2006}.
For a barrier parameter $\mu > 0$, IPM solves the smooth system of nonlinear equations
$F_\mu(w) = 0$ for $(s, \nu) > 0$ and $F_\mu(\cdot)$ defined as
\begin{equation}
  \label{eq:smoothkkt}
  F_\mu(w) = \begin{bmatrix}
     \nabla_x L(x, s, y, z, \nu) \\
     g(x) \\
     h(x) + s \\
     S \nu - \mu e_{m_i}\\
  \end{bmatrix} \; .
\end{equation}
We set $S = \text{diag}(s)$, $V = \text{diag}(\nu)$, and $e_{m_i} \in \mathbb{R}^{m_i}$ a
vector filled with $1$.
As we drive $\mu \to 0$, we recover the original KKT conditions~\eqref{eq:kkt}.

\subsection{Newton method}
The system of nonlinear equations~\eqref{eq:smoothkkt} is solved
using a Newton method. The primal-dual variable is updated as $w_{k+1} = w_k + \alpha d_k$, where
$d_k$ is a descent direction solution of the linear system
\begin{equation}
  \label{eq:newtonstep}
  \nabla_w F_\mu(w_k) d_k  = - F_\mu(w_k) \;.
\end{equation}
The step $\alpha$ is computed using a fraction-to-boundary rule, guaranteeing
that $(s, \nu) > 0$ throughout the iterations~\cite{nocedal_numerical_2006}.

\subsection{Augmented KKT system}
We note the local sensitivities $W_k = \nabla_{x x}^2 L(w_k)
\in \mathbb{R}^{n \times n}$,
$H_k = \nabla_x h(x_k) \in \mathbb{R}^{m_i \times n}$ and $G_k = \nabla_x g(x_k) \in \mathbb{R}^{m_e \times n}$.
The solution of the linear system~\eqref{eq:newtonstep}
translates to the \emph{augmented KKT system}:
\begin{equation}
  \label{eq:kkt:augmented}
  \overbrace{
  \begin{bmatrix}
    W_k  & 0 & G_k^\top & H_k^\top \\
    0 & D_s & 0& I \\
    G_k & 0 & 0 & 0 \\
    H_k & I & 0 & 0
\end{bmatrix}}^{K_{aug}}
  \begin{bmatrix}
    d_x \\
    d_s \\
    d_y \\
    d_z
  \end{bmatrix}
  = - \begin{bmatrix}
    r_1 \\ r_2 \\ r_3 \\ r_4
  \end{bmatrix} \; ,
\end{equation}
with the diagonal matrix $D_s = S_k^{-1} V_k$.
The right-hand-sides are given respectively by
$r_1 = \nabla f(x_k) + \nabla g(x_k)^\top y_k + \nabla h(x_k)^\top z_k + \mu X^{-1} e$,
$r_2 = z_k + \mu S^{-1} e$,
$r_3 = g(x_k)$,
$r_4 = h(x_k) + s_k$.

The system~\eqref{eq:kkt:augmented} is sparse, symmetric and exhibits
a saddle-point structure.  Most nonlinear optimization solvers solve
the system~\eqref{eq:kkt:augmented} using a sparse LBL factorization~\cite{duff1983multifrontal}.
Using \cite[Theorem 3.4]{benzi2005numerical},
the system \eqref{eq:kkt:augmented} is invertible if the Jacobian $J_k = \begin{bmatrix}
  G_k & 0 \\ H_k & I \end{bmatrix}$ is full row-rank and
\begin{equation}
  \label{eq:invertcond}
  \text{null}\left(
    \begin{bmatrix} W_k & 0 \\ 0 & D_s \end{bmatrix}
  \right)
  \cap
  \text{null}\left(
  \begin{bmatrix} G_k & 0 \\ H_k & I \end{bmatrix}
  \right)
  = \{0 \} \; .
\end{equation}
To ensure \eqref{eq:invertcond} holds, the solver checks the inertia
$\text{In}(K_{aug})$ (the tuple $(n_+, n_0, n_-)$ encoding respectively
the number of positive, null and negative eigenvalues in $K_{aug}$). If
\begin{equation}
  \label{eq:inertia}
  \text{In}(K_{aug}) = (n + m_i, 0, m_i + m_e) \; ,
\end{equation}
then the system $K_{aug}$ is invertible and
the solution of the system~\eqref{eq:kkt:augmented}
is a descent direction. Otherwise,
the solver regularizes~\eqref{eq:kkt:augmented} using two parameters $(\delta_x, \delta_c) > 0$
and solves
\begin{equation}
  \begin{bmatrix}
    W_k + \delta_x I  & 0 & G_k^\top & H_k^\top \\
    0 & D_s + \delta_x I  & 0& I \\
    G_k & 0 & -\delta_c I & 0 \\
    H_k & I & 0 & -\delta_c I
\end{bmatrix} \; .
\end{equation}
The parameter $(\delta_x, \delta_c)$ are computed so as the regularized
system satisfies~\eqref{eq:inertia}.
There exists inertia-free variants for IPM~\cite{chiang2016inertia}, but experimentally,
inertia-based method are known to converge in fewer iterations.

\subsection{IPM and optimal control}
IPM is a standard method to solve MPC and optimal control
problems~\cite{rao1998application}.
In particular, if \eqref{eq:problem} encodes a problem with a dynamic structure,
the Hessian $W_k$ and the Jacobian $H_k$ are block diagonal,
the Jacobian $G_k$ playing the role of the coupling matrix.

There exist interesting refinements of the IPM method for problems
with a dynamic structure~\cite{frison2020hpipm}. Notably:
\begin{itemize}
  \item The primal regularization $\delta_x$ can be computed
    recursively using dynamic programming, in a way that
    keeps the primal solution of the system \eqref{eq:kkt:augmented} intact
    \cite{verschueren2017sparsity}.
  \item Similarly, the dual regularization $\delta_c$ can be refined
    to take into account the problem's structure, implicitly
    (inside the linear solver~\cite{wan2017structured}) or explicitly (using exact penalty~\cite{thierry2020l1}).
\end{itemize}

\section{CONDENSED KKT SYSTEM}
\label{sec:condensed}
Solving the augmented KKT system~\eqref{eq:kkt:augmented} is numerically
demanding, and is often the computational bottleneck in IPM.
Furthermore, the sparse LBL factorization is known to be non trivial
to parallelize, as it relies on extensive numerical pivoting
operations~\cite{swirydowicz2021linear}. Fortunately, the KKT system
\eqref{eq:kkt:augmented} can be reduced down to a positive
definite matrix, whose factorization can be computed efficiently
using a Cholesky factorization.

First, we exploit the structure of the system~\eqref{eq:kkt:augmented}
using a condensation step\footnote{
  Here, the condensation removes the blocks associated
  to the slacks and the inequalities in the KKT system.
  As such, it has a different meaning than the condensing
  procedure used in model predictive control, which
  eliminates the state variables at all time except at 0 to
  obtain a dense KKT system~\cite{frison2020hpipm}.
}.
The system~\eqref{eq:kkt:augmented} is reduced by removing the
blocks associated to the slack $d_s$ and to the inequality multiplier $d_z$.
We obtain the equivalent \emph{condensed KKT system},
\begin{equation}
  \label{eq:kkt:condensed}
  \overbrace{
  \begin{bmatrix}
    K_k & G_k^\top \\
    G_k & 0
\end{bmatrix}}^{K_{cond}}
  \begin{bmatrix}
    d_x \\ d_y
  \end{bmatrix}
  =
  -
  \begin{bmatrix}
    r_1 + H_k^\top(D_s r_4 - r_2) \\ r_3
  \end{bmatrix}
   \; ,
\end{equation}
with the \emph{condensed matrix} $K_k := W_k + \delta_x + H_k^\top D_s H_k$.
Using the solution of the system~\eqref{eq:kkt:condensed},
we recover the updates on the slacks and inequality multipliers with
$d_s = -r_4 - H_k d_x$ and $d_z = -r_2 - D_s d_s$.
We note that the condensed matrix $K_k$ retains the block structure
of the Hessian $W_k$ and Jacobian $H_k$.

Using Haynsworth's inertia additivity formula, we have the equivalence
\begin{multline}
  \text{In}(K_{aug}) = (n + m_i, 0, m_i + m_e) \iff \\
  \text{In}(K_{cond}) = (n, 0, m_e)  \;.
\end{multline}


Usually, the solver uses a sparse LBL factorization to solve the KKT system $K_{aug}$ or $K_{cond}$.
Here, we move one step further and reduce the condensed KKT system~\eqref{eq:kkt:condensed}
down to a (sparse) positive definite matrix, using either
Lifted-KKT~\cite{shin2023accelerating} or HyKKT~\cite{regev2023hykkt}.

\subsection{Solution 1: Lifted-KKT}
We observe in \eqref{eq:kkt:condensed} that without equality constraints,
we obtain a $n \times n$ system which is guaranteed to be positive definite
if the primal regularization parameter $\delta_x$ is chosen appropriately.
Hence, we relax the equality constraints in \eqref{eq:problem}
using a small relaxation parameter $\tau > 0$, and solve the relaxed problem
\begin{equation}
  \label{eq:problemrelaxation}
    \min_{x \in \mathbb{R}^n} \; f(x) \quad
    \text{s.t.} \quad - \tau \leq g(x) \leq \tau \;, ~ h(x) \leq 0  \; .
\end{equation}
The problem~\eqref{eq:problemrelaxation} has only inequality constraints.
After introducing slack variables, the condensed KKT system
\eqref{eq:kkt:condensed} reduces to
\begin{equation}
  \label{eq:liftedkkt}
    K_k d_x = - r_1 - H_k^\top(D_s r_4 - r_2) \; .
\end{equation}
Using inertia correction method, the parameter $\delta_x$
is set to a value high enough to render the matrix $K_k$
positive definite. As a result, it can be factorized
efficiently using a sparse Cholesky method.

\subsection{Solution 2: HyKKT}
A substitute method is to exploit directly the structure of
the condensed KKT system~\eqref{eq:kkt:condensed}, without reformulating
the initial problem~\eqref{eq:problem}. To do so, we observe that if $K_k$ were
positive definite, the solution of the system~\eqref{eq:kkt:condensed} can be evaluated
using the Schur complement $G_k K_k^{-1} G_k^\top$.
Unfortunately the original problem \eqref{eq:problem} is nonconvex:
we have to convexify it using an Augmented
Lagrangian technique. For $\gamma > 0$, we note the KKT system~\eqref{eq:kkt:condensed}
is equivalent to
\begin{equation}
  \label{eq:hykkt}
  \begin{bmatrix}
    K_k + \gamma G_k^\top G_k & G_k^\top \\
    G_k & 0
  \end{bmatrix}
  \begin{bmatrix}
    d_x \\ d_y
  \end{bmatrix}
  =
  - \begin{bmatrix}
    r_\gamma \\
    r_3
  \end{bmatrix} \; .
\end{equation}
with $r_\gamma := r_1 + H_k^\top (D_s r_4 - r_2) + \gamma G_k^\top r_3$.

We note $Z$ a basis of the null-space of the Jacobian $G_k$.
We know that if $Z^\top K_k Z$ is positive definite and $G_k$ is full row-rank,
then there exists a threshold value $\underline{\gamma}$ such that
for all $\gamma \geq \underline{\gamma}$, $K_\gamma := K_k + \gamma G_k^\top G_k$
is positive definite~\cite{debreu1952definite}. This fact is exploited
in the Golub and Greif method~\cite{golub2003solving}, which has been
recently revisited in \cite{regev2023hykkt}. If $K_\gamma$
is positive definite, we can solve the system~\eqref{eq:hykkt}
using a Schur-complement method, by computing the
dual descent direction $d_y$ as solution of
\begin{equation}
  \label{eq:schurcomp}
  (G_k K_\gamma^{-1} G_k^\top) d_y =  r_3 - G_k K_\gamma^{-1} r_\gamma \;.
\end{equation}
Then, we recover the primal descent direction as
$K_\gamma d_x = r_\gamma - G_k^\top d_y$.
The method is tractable for two main reasons.
First, the matrix $K_\gamma$ is positive definite, meaning
it can be factorized efficiently without numerical pivoting.
Second, the Schur complement system \eqref{eq:schurcomp} can be solved
using a conjugate gradient algorithm converging in only a few iterations. In
fact, the eigenvalues of $S_\gamma := G_k K_\gamma^{-1} G_k^\top$ converge to
$\frac{1}{\gamma}$ as $\gamma \to +\infty$~\cite{regev2023hykkt}, implying
that the conditioning of $S_\gamma$ converges to $1$ (a setting particularly favorable for
iterative methods).

However, the conditioning of $K_\gamma$ increases with the parameter
$\gamma$ and the values in the diagonal matrix $D_s$.
Fortunately, the impact of the ill-conditioning remains limited in IPM,
and we can recover accurate solution when
solving the system~\eqref{eq:schurcomp} (see e.g. \cite{wright1998ill}).

\section{IMPLEMENTATION}
We have implemented the algorithm in the Julia
language, using the modeler ExaModels and the interior-point solver MadNLP~\cite{shin2023accelerating}.
Except for a few exceptions, all the array data is exclusively resident on the device memory, and the
algorithm has been designed to run fully on the GPU to avoid expensive
data transfers between the host and the device.

\subsection{Evaluation of the model with ExaModels}
Despite being nonlinear, the problem~\eqref{eq:problem}
usually has a highly repetitive structure that eases its
evaluation. CasADi allows for fast evaluation of problems with dynamic
structures~\cite{andersson2019casadi} but is not compatible with
GPU. Instead, we use the modeler ExaModels~\cite{shin2023accelerating}, which
detects the repeated patterns inside a nonlinear program
to evaluate them in a vectorized fashion using SIMD parallelism.
Using the multiple dispatch feature of Julia, ExaModels generates
highly efficient derivative computation code, compiled
for each computational pattern found in the model. Derivative evaluation
is implemented via array and kernel programming in the Julia Language,
using the wrapper CUDA.jl to dispatch the evaluation on the GPU.

\subsection{Hybrid linear solver}
The two KKT solvers introduced in \S\ref{sec:condensed},
Lifted-KKT and HyKKT, both require a sparse Cholesky factorization
and an iterative routine. We use the sparse Cholesky
solver cuDSS~\cite{CUDSS}, recently released by NVIDIA.
The most expensive operation in cuDSS is the computation
of the symbolic factorization.
However, it is important to note that the symbolic factorization can
be only computed once and refactorized efficiently
if the matrix's sparsity pattern remains the same.
This setting is particularly favorable for IPM,
as the sparsity pattern of the condensed matrix $K_k$
is fixed throughout the iterations.
Furthermore, within MPC framework, the symbolic factorization can be
performed offline, and thus, does not affect the online computation time.

Lifted-KKT factorizes the matrix $K_k$ using
cuDSS, and uses the resulting factor to solve the
linear system~\eqref{eq:liftedkkt} using a backsolve. The matrix
$K_k$ becomes increasingly ill-conditioned as
we approach the solution. Hence, it
has to be refined afterwards using an iterative refinement
algorithm, to increase the accuracy of the descent direction.
We use Richardson iterations in the iterative refinement, as in \cite{wachter2006implementation}.
Lifted-KKT sets the relaxation parameter to $\tau = 10^{-6}$.

HyKKT also factorizes the matrix $K_\gamma$ with cuDSS.
The resulting factor is used afterwards to
evaluate the residual $r_\gamma$ and solve
the Schur complement system \eqref{eq:schurcomp}
using a conjugate gradient algorithm (we only evaluate
matrix-vector products $S_\gamma x$ to avoid computing
the full matrix $S_\gamma$).
HyKKT leverages the conjugate gradient algorithm
implemented in the GPU-accelerated library Krylov.jl~\cite{montoison2023krylov}
to solve the system \eqref{eq:schurcomp} entirely on the GPU.
As the conditioning of the Schur complement $S_\gamma$ improves
with the parameter $\gamma$, the CG method
converges in less than 10 iterations on average, and
does not require a preconditioner.
In practice, we set $\gamma = 10^7$.

\subsection{GPU-accelerated interior-point solver}
The two KKT solvers Lifted-KKT and HyKKT have been
implemented inside MadNLP~\cite{shin2023accelerating},
a nonlinear solver implementing the filter line-search interior-point
method~\cite{wachter2006implementation} in pure Julia.
MadNLP builds upon the library CUDA.jl
to dispatch the operations seamlessly on the GPU, and leverages
the libraries cuSPARSE and cuBLAS for basic linear algebra operations.
The assembling of the two matrices $K_k$ and $K_\gamma$ occurs
entirely on the GPU, using a custom GPU kernel.

\section{NUMERICAL RESULTS}
We assess the performance of Lifted-KKT and HyKKT on the GPU, by comparing
them with the performance we obtain with the HSL solvers ma27 and ma57 on the CPU.

\subsection{Optimization of a distillation column}
As a test case, we use the classical distillation column instance
from \cite{cervantes1998large}, here implemented with ExaModels.
  The size of the discretization grid is parameterized by $N$.
  The distillation column has 32 trays.
  The vapor-liquid equilibrium equations are encoded with
  algebraic equations, the evolution of the material balances on each tray
  with differential equations.
  We note
  $\alpha=1.6$ the constant relative volatility,
  $L_t$ the liquid flow rate in the rectification section,
  $S_t$ the liquid flow rate in the stripping section,
  $x_{n,t}$ the liquid-phase mole fraction at tray $n$ and time $t$,
  $y_{n,t}$ the vapor-phase mole fraction,
  $u_t$ the reflux ratio,
  $D=0.2$ the constant distillate flow rate,
  $F=0.4$ the constant feed flow rate.
  The objective weights are $\gamma = 1000$ and $\rho=1$.
  The problem minimizes the quadratic deviation from the setpoints
  $(\bar{x}_1, \bar{u})$:
  \begin{equation*}
      \min  \sum_{t=1}^N \big( \gamma (x_{1,t} - \bar{x}_1)^2 + \rho (u_t - \bar{u})^2\big) \\
  \end{equation*}
  subject to, for all $t=1, \cdots, N$ and
  for a fixed time-step $\Delta t := 10 / N$,
  \begin{equation*}
    \begin{aligned}
        & L_t = u_t D \,, \; V_t = L_t + D \, , \; S_t = F + L_t \,, \\
        & y_{n,t} = \frac{\alpha x_{n,t}}{1 + (\alpha - 1) x_{n,t}} \,,
         \quad \forall n \in \{1,..,32\} \; , \\
        & \dot{x}_{1,t} = \frac{1}{M_1} V_t (y_{2,t} - x_{1,t}) \\
        & \dot{x}_{n,t} = \frac{1}{M_n} \big( L_t (x_{n-1,t} - x_{n,t}) -V_t (y_{n,t} - y_{n+1,t}) \big) \\
        & \quad \forall n \in \{2,..,16\} \; , \\
        & \dot{x}_{17,t} = \frac{1}{M_{17}} \big(F x_f + L_t x_{16,t} - S_t x_{17,t} - V_t (y_{17,t} - y_{18,t}) \big) \\
        & \dot{x}_{n,t} = \frac{1}{M_n} \big( S_t (x_{n-1,t} - x_{n,t}) -V_t (y_{n,t} - y_{n+1,t}) \big) \\
        & \quad \forall n \in \{18,..,31\} \; , \\
        & \dot{x}_{32,t} = \frac{1}{M_{32}} \big( S_t (x_{31,t} - (F - D_t) x_{32,t} - V_t y_{32,t} \big) \\
        & \dot{x}_{n, t} = \frac{1}{\Delta t}\big( x_{n,t} - x_{n,t-1} \big) \;,
         \quad \forall n \in \{1,..,32\} \; , \\
        & x_{n,0} = \bar{x}_{n,0} \;, ~ 1 \leq u_t \leq 5 \, ,
    \end{aligned}
  \end{equation*}

\subsection{Results}
We use our local workstation, equipped with an AMD
Epyc 7443 (24-core, 3.1GHz) and a NVIDIA A30 (24GB of local memory).
\footnote{
A script to reproduce the results is available at
\url{https://github.com/exanauts/nlp-on-gpu-paper/tree/main/HybridKKT.jl/benchmarks/cdc}.
}

\subsubsection{Performance of the cuDSS solver}
We start by assessing the performance of the sparse Cholesky
solver cuDSS~\cite{CUDSS}. We use MadNLP to generate the condensed matrix $K_k$
for the distillation column instance, for different discretization sizes
$N$. We benchmark individually
(i) the time to perform the symbolic analysis,
(ii) the time to refactorize the matrix,
(iii) the time to compute the backsolve.
The results are displayed in Table~\ref{tab:cudss}.
We note that for the largest instance ($N = 50,000$) we have
more than 3.3M variables $n$ in the optimization problem, but
only 0.0002\% of nonzeros in the sparse matrix $K_k$.
We observe that cuDSS spends most of its time in the symbolic analysis:
once the symbolic factorization is computed, the refactorization and the backsolve
are surprisingly fast (less than 0.5 seconds to recompute the factorization
of the largest instance). In other words, the costs of the symbolic factorization
is amortized as soon as we use many refactorizations afterward.
Furthermore, during the online computation within control systems, the
symbolic factorization can simply be reused and does not affect the
online computation time.

\begin{table}[thpb]
  \centering
    \caption{Performance of the linear solver cuDSS
      \label{tab:cudss}}
  \resizebox{.48\textwidth}{!}{
  \begin{tabular}{lrrrrrr}
  \hline
  $N$ & $n$ & nnz & SYM (s) & FAC (s) & SOLVE (s)  \\
  \hline
  100 & 6,767 & 83,835 & 0.058 & 0.003 & 0.001 \\
  500 & 33,567 & 418,635 & 0.209 & 0.011 & 0.004 \\
  1,000 & 67,067 & 837,135 & 0.381 & 0.019 & 0.006 \\
  5,000 & 335,067 & 4,185,135 & 2.185 & 0.072 & 0.022 \\
  10,000 & 670,067 & 8,370,135 & 4.396 & 0.115 & 0.037 \\
  20,000 & 1,340,067 & 16,740,135 & 9.057 & 0.220 & 0.072 \\
  50,000 & 3,350,067 & 33,450,135 & 20.187 & 0.432 & 0.165 \\
  \hline
  \end{tabular}
  }
  \\[.25em]
  \footnotesize
  $^*$SYM, FAC, and SOLVE are the time spent in the symbolic
analysis, refactorization, and backsolve, respectively (all displayed in seconds).
\end{table}

\subsubsection{Performance of the MadNLP solver}
We analyze now the performance in the MadNLP solver, using both
the Lifted-KKT and HyKKT hybrid solvers together with cuDSS. As a baseline,
we give the time spent in the HSL solvers ma27 and ma57.
We set MadNLP tolerance to {\tt tol=1e-6}, and look at the time
spent to achieve convergence in IPM.
The results are displayed in Table~\ref{tab:madnlp}.
We observe that MadNLP converges in the same number of iterations
with HSL ma27, HSL ma57 and HyKKT, as these KKT solvers all solve the
same KKT system~\eqref{eq:kkt:augmented}. On its end, Lifted-KKT solves
the relaxed problem~\eqref{eq:problemrelaxation} and requires twice as much
iterations to achieve convergence. In concordance with Table~\ref{tab:cudss},
the pre-processing times (column {\tt init}) are significant for Lifted-KKT and HyKKT.
Computing the symbolic factorization with cuDSS is the bottleneck: in comparison,
ma27 is approximately twice as fast during the pre-processing. However, the
time spent in the symbolic factorization is amortized throughout the IPM iterations:
Lifted-KKT and HyKKT solve the problem respectively 26x and 18x faster than HSL ma27
on the largest instance ($N=50,000$).

\begin{figure}[thpb]
    \centering
    \includegraphics[width=.48\textwidth]{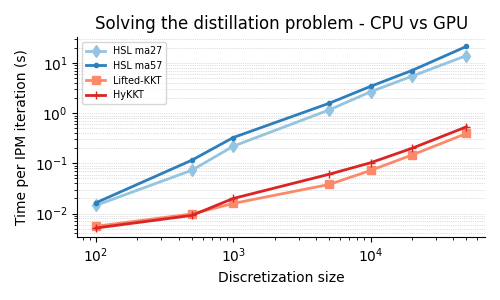}
    \caption{Time per IPM iteration (s), CPU versus GPU.}
    \label{fig:timeipmit}
\end{figure}

We display the time per IPM iteration in Figure~\ref{fig:timeipmit}.
We observe that Lifted-KKT is slightly faster than HyKKT, as the performance
of the later method depends on the total number of CG iterations required
to solve the system~\eqref{eq:schurcomp} at each IPM iteration.
The solver HSL ma27 is consistently better than ma57 on that particular instance,
as the problem is highly sparse.

\begin{table*}[t]
  \centering
  \caption{Performance comparison of MadNLP on CPU and GPU.
  }
  \resizebox{.9\textwidth}{!}{
  \begin{tabular}{l|rrr>{\bfseries}r|rrr>{\bfseries}r|rrr>{\bfseries}r|rrr >{\bfseries}r}
    \hline
  & \multicolumn{4}{c|}{HSL ma27} &
    \multicolumn{4}{c|}{HSL ma57} &
    \multicolumn{4}{c|}{Lifted-KKT} &
    \multicolumn{4}{c}{HyKKT} \\
    \hline
\hline
 $N$ & init & AD & linsolve & total & init & AD & linsolve & total & init & AD & linsolve & total & init & AD & linsolve & total \\
\hline
100 & 0.0 & 0.0 & 0.1 & 0.1 & 0.0 & 0.0 & 0.1 & 0.1 & 0.1 & 0.0 & 0.2 & 0.3 & 0.1 & 0.0 & 0.1 & 0.2 \\
500 & 0.1 & 0.0 & 0.6 & 0.7 & 0.0 & 0.0 & 1.0 & 1.0 & 0.2 & 0.0 & 0.3 & 0.6 & 0.2 & 0.0 & 0.1 & 0.3 \\
1,000 & 0.1 & 0.0 & 1.7 & 1.8 & 0.6 & 0.0 & 2.3 & 3.0 & 0.5 & 0.0 & 0.4 & 0.9 & 0.4 & 0.0 & 0.2 & 0.6 \\
5,000 & 0.6 & 0.1 & 8.8 & 9.5 & 3.6 & 0.1 & 13.2 & 16.9 & 2.3 & 0.0 & 0.5 & 2.8 & 2.3 & 0.0 & 0.4 & 2.7 \\
10,000 & 1.6 & 0.2 & 20.6 & 22.4 & 7.7 & 0.2 & 26.6 & 34.4 & 4.8 & 0.0 & 0.9 & 5.8 & 4.9 & 0.0 & 0.8 & 5.7 \\
20,000 & 4.6 & 0.4 & 40.3 & 45.3 & 17.4 & 0.4 & 59.0 & 76.8 & 10.2 & 0.1 & 2.0 & 12.3 & 10.6 & 0.1 & 2.0 & 12.6 \\
50,000 & 15.4 & 0.9 & 109.1 & 125.5 & 49.5 & 0.9 & 146.3 & 196.8 & 27.9 & 0.5 & 5.4 & 33.8 & 29.7 & 0.1 & 4.6 & 34.5 \\
\hline
\end{tabular}
}
  \label{tab:madnlp}
\end{table*}

\section{CONCLUSIONS AND FUTURE WORKS}

\subsection{Conclusions}
In this paper, we have presented two hybrid solvers to solve
the sparse KKT systems arising in IPM on the GPU: Lifted-KKT and HyKKT. We have implemented
the two hybrid solvers within MadNLP, using the newly released cuDSS linear solver
to compute the sparse Cholesky factorization on the GPU.
Our results on the distillation column instance show that once the symbolic factorization is computed,
both Lifted-KKT and HyKKT are significantly faster than HSL running on the CPU,
with a time per iteration reduced by a factor of 25 on the largest instance.
This setting is relevant for NMPC, as the symbolic factorization
can be computed only once and reused when we solve the problem
a second time with updated data.

\subsection{Future Works}
In this article, the dynamic structure
has not been exploited explicitly. However, it is known that
for dynamic nonlinear programs, the condensed matrix $K_k$
is block-banded and can be factorized efficiently
using a block-structured linear solver. Indeed, we can
interpret $K_k$ as a matrix encoding a linear-quadratic program
and solve the resulting problem in parallel using partitioned dynamic
programming~\cite{wright1991partitioned}.
We are planning to investigate this promising outlook
in future research.

\addtolength{\textheight}{-2.7cm}   



\begin{thebibliography}{10}

\bibitem{diehl2009efficient}
M.~Diehl, H.~J. Ferreau, and N.~Haverbeke, ``Efficient numerical methods for nonlinear {MPC} and moving horizon estimation,'' {\em Nonlinear model predictive control: towards new challenging applications}, pp.~391--417, 2009.

\bibitem{kirches2010efficient}
C.~Kirches, L.~Wirsching, S.~Sager, and H.~G. Bock, ``Efficient numerics for nonlinear model predictive control,'' in {\em Recent Advances in Optimization and its Applications in Engineering: The 14th Belgian-French-German Conference on Optimization}, pp.~339--357, Springer, 2010.

\bibitem{dunn1989efficient}
J.~C. Dunn and D.~P. Bertsekas, ``Efficient dynamic programming implementations of {N}ewton's method for unconstrained optimal control problems,'' {\em Journal of Optimization Theory and Applications}, vol.~63, no.~1, pp.~23--38, 1989.

\bibitem{zavala2008interior}
V.~M. Zavala, C.~D. Laird, and L.~T. Biegler, ``Interior-point decomposition approaches for parallel solution of large-scale nonlinear parameter estimation problems,'' {\em Chemical Engineering Science}, vol.~63, no.~19, pp.~4834--4845, 2008.

\bibitem{wright1991partitioned}
S.~J. Wright, ``Partitioned dynamic programming for optimal control,'' {\em SIAM Journal on Optimization}, vol.~1, no.~4, pp.~620--642, 1991.

\bibitem{rao1998application}
C.~V. Rao, S.~J. Wright, and J.~B. Rawlings, ``Application of interior-point methods to model predictive control,'' {\em Journal of Optimization Theory and Applications}, vol.~99, no.~3, pp.~723--757, 1998.

\bibitem{cervantes1998large}
A.~Cervantes and L.~T. Biegler, ``Large-scale {DAE} optimization using a simultaneous {NLP} formulation,'' {\em AIChE Journal}, vol.~44, no.~5, pp.~1038--1050, 1998.

\bibitem{andersson2019casadi}
J.~A. Andersson, J.~Gillis, G.~Horn, J.~B. Rawlings, and M.~Diehl, ``{CasADi}: a software framework for nonlinear optimization and optimal control,'' {\em Mathematical Programming Computation}, vol.~11, pp.~1--36, 2019.

\bibitem{verschueren2022acados}
R.~Verschueren, G.~Frison, D.~Kouzoupis, J.~Frey, N.~v. Duijkeren, A.~Zanelli, B.~Novoselnik, T.~Albin, R.~Quirynen, and M.~Diehl, ``acados—a modular open-source framework for fast embedded optimal control,'' {\em Mathematical Programming Computation}, vol.~14, no.~1, pp.~147--183, 2022.

\bibitem{ferreau2014qpoases}
H.~J. Ferreau, C.~Kirches, A.~Potschka, H.~G. Bock, and M.~Diehl, ``{qpOASES}: A parametric active-set algorithm for quadratic programming,'' {\em Mathematical Programming Computation}, vol.~6, pp.~327--363, 2014.

\bibitem{frasch2015parallel}
J.~V. Frasch, S.~Sager, and M.~Diehl, ``A parallel quadratic programming method for dynamic optimization problems,'' {\em Mathematical programming computation}, vol.~7, pp.~289--329, 2015.

\bibitem{frison2020hpipm}
G.~Frison and M.~Diehl, ``{HPIPM}: a high-performance quadratic programming framework for model predictive control,'' {\em IFAC-PapersOnLine}, vol.~53, no.~2, pp.~6563--6569, 2020.

\bibitem{frison2018blasfeo}
G.~Frison, D.~Kouzoupis, T.~Sartor, A.~Zanelli, and M.~Diehl, ``{BLASFEO}: Basic linear algebra subroutines for embedded optimization,'' {\em ACM Transactions on Mathematical Software (TOMS)}, vol.~44, no.~4, pp.~1--30, 2018.

\bibitem{kerrigan2015computer}
E.~C. Kerrigan, G.~A. Constantinides, A.~Suardi, A.~Picciau, and B.~Khusainov, ``Computer architectures to close the loop in real-time optimization,'' in {\em 2015 54th IEEE conference on decision and control (CDC)}, pp.~4597--4611, IEEE, 2015.

\bibitem{phung2017model}
D.-K. Phung, B.~H{\'e}riss{\'e}, J.~Marzat, and S.~Bertrand, ``Model predictive control for autonomous navigation using embedded graphics processing unit,'' {\em IFAC-PapersOnLine}, vol.~50, no.~1, pp.~11883--11888, 2017.

\bibitem{yu2017efficient}
L.~Yu, A.~Goldsmith, and S.~Di~Cairano, ``Efficient convex optimization on {GPUs} for embedded model predictive control,'' in {\em Proceedings of the General Purpose GPUs}, (New York, NY, USA), pp.~12--21, Association for Computing Machinery, 2017.

\bibitem{rathai2020gpu}
K.~M.~M. Rathai, M.~Alamir, and O.~Sename, ``{GPU} based stochastic parameterized {NMPC} scheme for control of semi-active suspension system for half car vehicle,'' {\em IFAC-PapersOnLine}, vol.~53, no.~2, pp.~14369--14374, 2020.

\bibitem{shin2023accelerating}
S.~Shin, M.~Anitescu, and F.~Pacaud, ``Accelerating optimal power flow with {GPUs}: {SIMD} abstraction of nonlinear programs and condensed-space interior-point methods,'' {\em Electric Power Systems Research}, vol.~236, p.~110651, 2024.

\bibitem{wachter2006implementation}
A.~W{\"a}chter and L.~T. Biegler, ``On the implementation of an interior-point filter line-search algorithm for large-scale nonlinear programming,'' {\em Mathematical Programming}, vol.~106, no.~1, pp.~25--57, 2006.

\bibitem{benzi2005numerical}
M.~Benzi, G.~H. Golub, and J.~Liesen, ``Numerical solution of saddle point problems,'' {\em Acta numerica}, vol.~14, pp.~1--137, 2005.

\bibitem{duff1983multifrontal}
I.~S. Duff and J.~K. Reid, ``The multifrontal solution of indefinite sparse symmetric linear,'' {\em ACM Transactions on Mathematical Software (TOMS)}, vol.~9, no.~3, pp.~302--325, 1983.

\bibitem{tasseff2019exploring}
B.~Tasseff, C.~Coffrin, A.~W{\"a}chter, and C.~Laird, ``Exploring benefits of linear solver parallelism on modern nonlinear optimization applications,'' {\em arXiv preprint arXiv:1909.08104}, 2019.

\bibitem{swirydowicz2021linear}
K.~{\'S}wirydowicz, E.~Darve, W.~Jones, J.~Maack, S.~Regev, M.~A. Saunders, S.~J. Thomas, and S.~Pele{\v{s}}, ``Linear solvers for power grid optimization problems: a review of {GPU}-accelerated linear solvers,'' {\em Parallel Computing}, p.~102870, 2021.

\bibitem{cole2023exploiting}
D.~Cole, S.~Shin, F.~Pacaud, V.~M. Zavala, and M.~Anitescu, ``Exploiting {GPU/SIMD} architectures for solving linear-quadratic {MPC} problems,'' in {\em 2023 American Control Conference (ACC)}, pp.~3995--4000, IEEE, 2023.

\bibitem{pacaud2022condensed}
F.~Pacaud, S.~Shin, M.~Schanen, D.~A. Maldonado, and M.~Anitescu, ``Accelerating condensed interior-point methods on {SIMD/GPU} architectures,'' {\em Journal of Optimization Theory and Applications}, pp.~1--20, 2023.

\bibitem{regev2023hykkt}
S.~Regev, N.-Y. Chiang, E.~Darve, C.~G. Petra, M.~A. Saunders, K.~{\'S}wirydowicz, and S.~Pele{\v{s}}, ``{HyKKT}: a hybrid direct-iterative method for solving {KKT} linear systems,'' {\em Optimization Methods and Software}, vol.~38, no.~2, pp.~332--355, 2023.

\bibitem{golub2003solving}
G.~H. Golub and C.~Greif, ``On solving block-structured indefinite linear systems,'' {\em SIAM Journal on Scientific Computing}, vol.~24, no.~6, pp.~2076--2092, 2003.

\bibitem{nocedal_numerical_2006}
J.~Nocedal and S.~J. Wright, {\em Numerical optimization}.
\newblock Springer series in operations research, New York: Springer, 2nd~ed., 2006.

\bibitem{chiang2016inertia}
N.-Y. Chiang and V.~M. Zavala, ``An inertia-free filter line-search algorithm for large-scale nonlinear programming,'' {\em Computational Optimization and Applications}, vol.~64, pp.~327--354, 2016.

\bibitem{verschueren2017sparsity}
R.~Verschueren, M.~Zanon, R.~Quirynen, and M.~Diehl, ``A sparsity preserving convexification procedure for indefinite quadratic programs arising in direct optimal control,'' {\em SIAM Journal on Optimization}, vol.~27, no.~3, pp.~2085--2109, 2017.

\bibitem{wan2017structured}
W.~Wan and L.~T. Biegler, ``Structured regularization for barrier {NLP} solvers,'' {\em Computational Optimization and Applications}, vol.~66, pp.~401--424, 2017.

\bibitem{thierry2020l1}
D.~Thierry and L.~Biegler, ``The $\ell$1—exact penalty-barrier phase for degenerate nonlinear programming problems in {Ipopt},'' {\em IFAC-PapersOnLine}, vol.~53, no.~2, pp.~6496--6501, 2020.

\bibitem{debreu1952definite}
G.~Debreu, ``Definite and semidefinite quadratic forms,'' {\em Econometrica: Journal of the Econometric Society}, pp.~295--300, 1952.

\bibitem{wright1998ill}
M.~H. Wright, ``Ill-conditioning and computational error in interior methods for nonlinear programming,'' {\em SIAM Journal on Optimization}, vol.~9, no.~1, pp.~84--111, 1998.

\bibitem{CUDSS}
{NVIDIA cuDSS documentation}, ``{NVIDIA cuDSS (Preview): A high-performance CUDA Library for Direct Sparse Solvers}.''
\newblock \url{https://docs.nvidia.com/cuda/cudss/index.html}.

\bibitem{montoison2023krylov}
A.~Montoison and D.~Orban, ``Krylov. jl: A julia basket of hand-picked krylov methods,'' {\em Journal of Open Source Software}, vol.~8, no.~89, p.~5187, 2023.

\end{thebibliography}
\end{document}